 \documentclass[12pt]{article}

\usepackage{amssymb}
\def\vep{\varepsilon}

\def\reff#1{(\ref{#1})}
\def\gbt{{t\gamma\beta }}

\def\E{{\bf E}}
\def\P{{\bf P}}
\def\R{{\mathbb R}}
\def\Z{{\mathbb Z}}
\def\N{{\mathbb N}}
\def\X{{\bf X}}
\def\1{{\mathbb 1}}

\def\ga{{\gamma}}

\def\proof{\noindent{\bf Proof. }}

\def\B{{\bf B}}

\def\one{{\bf 1}}

\begin{document}

\title{Convergence to the maximal invariant measure\\
 for a zero-range process with random rates.}

\author{E.D. Andjel\kern -0.5pt
\renewcommand{\thefootnote}{\alph{footnote}}\footnotemark%
\ \ P.A. Ferrari\kern -2pt \addtocounter{footnote}{4}
\renewcommand{\thefootnote}{\alph{footnote}}\footnotemark%
\ \ H. Guiol\kern -0.5pt
\renewcommand{\thefootnote}{\alph{footnote}}\footnotemark%
\ \ and\ C. Landim\kern -2pt
\addtocounter{footnote}{4}\renewcommand{\thefootnote}{\alph{footnote}}
\footnotemark}

\maketitle
\renewcommand{\thefootnote}{\alph{footnote}}
\addtocounter{footnote}{1} \footnotetext{LATP-CMI, 36 Rue
Joliot-Curie, 13013 Marseille, France.}
\renewcommand{\thefootnote}{\alph{footnote}}
\addtocounter{footnote}{5} \footnotetext{IME-USP, P.B. 66281,
05315-970 S\~ao Paulo, SP, Brasil.}
\renewcommand{\thefootnote}{\alph{footnote}}
\addtocounter{footnote}{1} \footnotetext{IMECC-UNICAMP, P.B. 6065,
13053-970, Campinas, SP, Brasil.}
\renewcommand{\thefootnote}{\alph{footnote}}
\addtocounter{footnote}{5} \footnotetext{IMPA, Estrada Dona
Castorina 110, Jardim Bot\^anico, Rio de Janeiro, Brasil and CNRS
UPRES-A 6085, Universit\'e de Rouen, BP 118, 76821 Monts Saint
Aignan Cedex, France.}

\newcommand{\carn}{\hfill\rule{0.25cm}{0.25cm}}

\newtheorem{theorem}{Theorem}[section]
\newtheorem{lemma}[theorem]{Lemma}
\newtheorem{proposition}[theorem]{Proposition}
\newtheorem{corollary}[theorem]{Corollary}
\newtheorem{conjecture}[theorem]{Conjecture}
\newtheorem{definition}[theorem]{Definition}
\newtheorem{remark}[theorem]{Remark}

\abstract {We consider a one-dimensional totally asymmetric
  nearest-neighbor zero-range process with site-dependent jump-rates
  ---an \emph{environment}. For each environment $p$ we prove that the
  set of all invariant measures is the convex hull of a set of product
  measures with geometric marginals. As a consequence we show that
  for environments $p$ satisfying certain asymptotic
  property, there are no invariant measures concentrating on
  configurations with density bigger than $\rho^*(p)$, a
  critical value. If $\rho^*(p)$ is finite we say that there is
  phase-transition on the density. In this case we prove that if the
  initial configuration has asymptotic density strictly above
  $\rho^*(p)$, then the process converges to the
  maximal invariant measure.
   }

\section{Introduction}

The interest on the behavior of interacting particle systems in
random environment has grown recently: Benjamini, Ferrari and
Landim (1996), Evans (1996) and Ferrari and Krug (1996) observed
the existence of phase transition in these models; Benjamini,
Ferrari and Landim (1996), Krug and Sepp\"al\"ainen (1999) and
Koukkous (1999) investigated the hydrodynamic behavior of
conservative processes in random environments; Landim (1996) and
Bahadoran (1998) considered the same problem for non-homogeneous
asymmetric attractive processes; Gielis, Koukkous and Landim
(1998) deduced the equilibrium fluctuations of a symmetric zero
range process in a random environment.

In this article we consider a one-dimensional, totally asymmetric,
nearest-neighbor zero-range process in a non-homogeneous
environment. The evolution can be informally described as follows.
Fix $c\in(0,1)$ and provide each site $x$ of $\Z$ with a rate
function $p_x\in[c,1]$. If there is at least one particle at some
site $x$, one of these particles jumps to $x+1$ at rate $p_x$. A
rate configuration $p=(p_x:x\in\Z)$ is called an \emph
{environment} and a measure $m$ on the set of possible
environments a \emph{random environment.}

Benjamini, Ferrari and Landim (1996) and Evans (1996) for an
asymmetric exclusion process with rates associated to the
particles
---which is isomorphic to a zero range process with rates associated
to the sites--- and Ferrari and Krug (1996) for the model
considered here, proved the existence of a phase transition in the
density. More precisely, they proved that, under certain
conditions on the distribution $m$, specified in Theorem
\ref{2.4}, there exists a finite critical value $\rho^*$ such that
for $m$-almost-all $p$ there are no product invariant measures for
the process with rates $p$ concentrating on configurations with
asymptotic density bigger than $\rho^*$ and that there are product
invariant measures concentrating on configurations with asymptotic
density smaller than or equal to $\rho^*$. (The density of a
configuration is essentially the average number of particles per
site and is defined in \reff{dens} below).

Our first result is that the set of extremal invariant measures
for the process with fixed environment $p=(p_x:x\in\Z)$ is the set
$\{\nu_{p,v}: v< p_x, \forall x\}$, where $\nu_{p,v}$ is the
product measure on $\N^\Z$ with marginals
\begin{equation}
\label{889} \nu_{p,v}\{\xi : \, \xi(x)=k\} =\Bigl(\frac v{
p_x}\Bigr)^k \Bigl(1-\frac v{ p_x}\Bigr)\; .
\end{equation}
The above result does not surprise specialists in queuing theory.
In fact we are dealing with an infinite series of M/M/1 queues
with service rate $p_x$ at queue $x$. The value $v$ can be
interpreted as the arrival rate at ``queue'' $-\infty$. Since
Burke's theorem (see Kelly (1979) or Theorem 7.1 in Ferrari (1992)
for instance) guarantees that in equilibrium the departure process
of a M/M/1 queue is the same as the arrival process (both Poisson
of rate $v$), there is an invariant measure for each arrival rate
$v$ strictly smaller than all service rates.

Assume $c=\inf_x p_x$ and that the following limits exist. For
$v<c$,
\begin{equation}
  \label{cc1}
  R(p,v)\,:=\, \lim_{n\to\infty}{1\over n}\sum_{x=-n+1}^0 \int
  \nu_{p,v}(d\xi)\, \xi(x)\, = \,\lim_{n\to\infty}{1\over n}
  \sum_{x=-n+1}^0 \frac {v}{p_x -v}.
\end{equation}
We interpret $R(p,v)$ as the global expected left density per site
of the configurations distributed according to $\nu_{p,v}$. A
consequence of the existence of the limits, as we will explain
later, is that for all $v<c$, $\nu_{p,v}$ concentrates in
configurations with asymptotic left density $R(p,v)$:
\begin{equation}
  \label{cc5}
  \nu_{p,v}\Bigl(\lim_{n\to\infty}{1\over n}\sum_{x=-n+1}^0 \xi(x) =
R(p,v)\Bigr)\,=\,1\; .
\end{equation}
It is easy to prove that $R(p,v)$ is a strictly convex increasing
function of $v$, hence the limit
\begin{equation}
  \label{cc6}
  \rho^*(p):=\lim_{v\to c} R(p,v)
\end{equation}
is well defined (but may be infinite). In the sequel we assume
$\rho^*(p)<\infty$. We do not assume the existence of the limit in
\reff{cc1} for $v=c$, nor the $\nu_{p,c}$ almost sure convergence
of the density.

Our second and main result states that under the condition
$\rho^*(p)<\infty$, initial measures concentrating on
configurations with asymptotic left density strictly bigger than
$\rho^*(p)$ converge towards the maximal invariant measure
$\nu_{p,c}$. We do not know in general if this measure
concentrates on configurations with density. But if the limit
$R(p,c)$ of \reff{cc1} exists, equals $\rho^*(p)$ and is finite
our result says that the process starts with global density
strictly above $\rho^*(p)$ and converges to a measure with density
$\rho^*(p)$. This behavior is remarkable as the process is
\emph{conservative}, i.e. the total number of particles is
conserved, but in the above limit ``looses mass''.  Informally
speaking, what happens is that many clients remain trapped in far
away slow servers. More precisely, denoting by $S_p(t)$ the
semigroup of the process, we first show that for any initial
measure $\nu$, all weak limits of the sequence $\{ \nu S_p(t),\,
t\ge 0\}$ are dominated, in the natural partial order, by
$\nu_{p,c}$.  We then show that if $\nu$ is a measure concentrated
on configurations with asymptotic left density strictly greater
than $\rho^*(p)$, all weak limits of $\nu S_p(t)$ dominate
$\nu_{p,c}$. Surprisingly enough, the proof of the second
statement is much more demanding than the proof of the first one.

It follows from the two previous results that the domain of
attraction of $\nu_{p,c}$ includes all measures with asymptotic
density strictly above $\rho^*(p)$. It remains an open question to
describe the domain of attraction of a product invariant measure
$\nu_{p,v}$ for $0< v<c$ or to show the convergence to $\nu_{p,c}$
of initial measures with asymptotic density $\rho^*(p)$.

Our results hold $m$-a.s.\ for measures $m$ concentrating on
environments satisfying \reff{cc1}.

The paper is organized as follows: in Section 2 we introduce the
notation and state the main results. In Section 3 we characterize
the set of invariant measures and show that the maximal invariant
measure dominates all the weak limits of the process. In Section 4
we obtain the asymptotic velocity of a second class particle for
the zero-range process in a non homogeneous environment and use
this result to prove the main theorem.

Many of our results are based on standard coupling arguments. We
assume the reader familiar with this technique described in
Section 1 of Chapter 2 of Liggett (1985).

\section{Notation and Results}

Fix $0<c\leq 1$ and consider a sequence $(p_x)_{x\in \Z}$ taking
values in $[c,1]$ such that $c=\inf_x p_x$. We consider a totally
asymmetric zero-range process in the environment $p$. This is a
Markov process that can be informally described as follows. We
initially distribute particles on the lattice $\Z$ . If there is
at least one particle at some site $x$, then at rate $ p_x$ one of
them jumps to site $x+1$. To construct a Markov process $\eta_t$
on $\X=\N^{\Z}$ corresponding to the above description, let
$N_x(t)\ ( x\in \Z)$ be a collection of independent Poisson
processes such that for all $x\in \Z,\ \E(N_x(t))=p_xt$. The
evolution of $\eta_t$ is now given by the following rule: if the
Poisson process $N_x(. )$ jumps at time $t$ and $\eta_{t-}(x)>0$
then one particle is moved from $x$ to $x+1$ at that time. To see
that the process is well defined by this rule, just note that in
any time interval $[0,t]$ for any $ x$ there exists with
probability $1$ a $y<x$ such that $N_y(t)=0$.  Hence the value of
$\eta_t(x)$ depends only on the initial configuration and on a
finite number of jumps.

The generator $L_p$ of this process, defined by $L_pf(\eta) =
d\E[f(\eta_t) \,|\,\eta_0=\eta]/dt\Big|_{t=0}$, acts on cylinder
functions $f$ as follows:
\begin{equation}
 \label{a3}
(L_p f)(\eta)=\sum_{x\in{\Z}}\,p_x\,\one\{\eta (x)>0\}
\left[f(\eta^x)-f(\eta)\right]\;.
\end{equation}
In the above formula $\eta^x=\eta-\mathfrak d_{x}+\mathfrak
d_{x+1}$, where $\mathfrak d_y$ stands for a configuration with
just one particle at $y$ and addition of configurations is
performed componentwise.

We denote by $\{S_p(t),\, t\ge 0\}$ the semigroup associated to
the generator $L_p$, i.e.\/ $S_p(t)f(\eta) =
\E[f(\eta_t)\,|\,\eta_0=\eta]$ and by ${\cal
  I}_p$ the set of invariant measures of $\eta_t$ (the Markov process with
generator $L_p$).  Let $v$ be a real number such that $ 0<v<p_x$
for all $x$. Then a standard calculation (first observed by
Jackson (1957) for the finite case) shows that the product measure
$\nu_{p,v}$ with marginals given by (\ref{889})
is an invariant measure for the process. Benjamini, Ferrari and
Landim (1996) raised the question of whether or not there exist
invariant measures which are not convex combinations of the
$\nu_{p,v}$'s.  In Section 3 we prove the following theorem which,
combined with Theorem 12.2 in Dynkin (1978) (which states that the
set of extremal invariant measures is the convex hull of the set
of invariant measures), gives a negative answer to that question.
In its statement we denote by $({\cal I}_p)_e$ the set of extremal
invariant measures for the process.
\begin{theorem}
\label{2.1} Let $p$ be an arbitrary environment then
 \[
 ({\cal I}_p)_e=\{\nu_{p,v}\ :\ v<p_x,\ \forall x\in \Z \}.
 \]
\end{theorem}

In this theorem the range of the parameter $v$ may be either
$[0,c)$ or $[0,c]$ - when $p_x=c$ for some $x$ or $p_x>c$ for all
$x$, respectively. In the first case to prove the theorem we only
need to follow the proof of Theorem 1.11 in Andjel (1982), but in
the second case a complementary argument is needed.  In both cases
the proof relies on the standard partial order for probability
measures on $\X$. To define it, first say that $\eta \leq \xi$ if
$\eta (x) \leq \xi (x)$ for all $x\in \Z$. Then say that a real
valued function $f$ defined on $\X$ is increasing if $ \eta \leq
\xi$ implies that $f(\eta) \leq f(\xi)$. Finally if $\mu$ and $\nu
$ are two probability measures on $\X$, say that $\mu \leq \nu$ if
$\int fd\mu \leq \int fd\nu$ for all bounded increasing cylinder
functions $f$. In this case we say that $\nu$ \emph{dominates}
$\mu$.  The complementary argument alluded above depends on the
following proposition:

\begin{proposition}
\label{2.2} Assume that $p$ is an environment such that
\begin{equation}
  \label{ccc}
p_x>c \mbox{ for all }x\in\Z\mbox{ and }\liminf_{x\to -\infty}\
p_x=c\,,
\end{equation}
and let $\nu$ be an arbitrary probability measure on $\X$. Then
the set of measures $\{\nu S_p(t)\; :\; t>0\}$ is tight and its
weak limits as $t$ goes to infinity are bounded above by
$\nu_{p,c}$.
\end{proposition}

An immediate corollary of Proposition \ref{2.2} is that under
\reff{ccc} all invariant measures are dominated by~$\nu_{p,c}$.

To state our main result let $\eta$ be an element of $\X$ and
consider
\begin{eqnarray}
  \label{dens}
  \underline D(\eta )&=& \liminf_{n\to\infty} \frac
  {1}{n}\sum_{x=-n+1}^0 \eta (x)\;,\nonumber\\
  \overline D(\eta )&=& \limsup_{n\to\infty} \frac
  {1}{n}\sum_{x=-n+1}^0 \eta (x)\;,\nonumber
\end{eqnarray}
the \emph{lower}, respectively \emph{upper asymptotic left
density} of $\eta$. If both limits are equal to $\alpha$ we say
that $\eta$ has \emph{left
  density} $\alpha$ and write $D(\eta)= \alpha$.

Assume that $p$ is an environment for which the limits defined in
\reff{cc1} exist. Then, by Kolmogorov's law of large numbers (see
{\sl e.g.} Shiryayev (1984), Theorem 2 p. 364) $\nu_{p,v}$
concentrates on configurations with left density $R(p,v)$:
\begin{equation}
  \label{cc2}
  \nu_{p,v}\{\eta\in\X: D(\eta) = R(p,v)\}= 1
\end{equation}
for all $v<c$.

The values assumed by $R(p,v)$ for $v<c$ are crucial for the
characterization of the set of invariant measures for the process
with rates $p$. If $\lim_{v\to c}R(p,v)=\infty$, then the range of
allowed densities is $[0,\infty)$ or $[0,\infty]$. The first case
occurs when $p_x= c $ for some~$x$. In this case $\nu_{p,v}$ is
defined for any $v<c$, but not for $v=c$.  Moreover, since $R(p,
\cdot)$ is continuous and increases to $\infty$ as $v\to c$, then
for all $\rho\in[0,\infty)$ there exists $v=v(p,\rho)$ such that
$\nu_{p,v}\{\eta\in\X: D(\eta) = \rho\}= 1$. The second case
occurs when $p_x>c$ for all $x$. In this case $\nu_{p,c}$ is well
defined and concentrates on configurations with infinite
asymptotic left density, and for any $\rho\in[0,\infty]$ there
exists $v=v(p,\rho)$ such that $\nu_{p,v}\{\eta\in\X: D(\eta) =
\rho\}= 1$.

If $\lim_{v\to c}R(p,v)= \rho^*(p)<\infty$ and $p_x>c$ for all
$x$, the measure $\nu_{p,c}$ is well defined and Theorem \ref{2.1}
tells us that there are no invariant measures  bigger than
$\nu_{p,c}$. Our next theorem describes what happens in this case
when one starts with a density strictly bigger than $\rho^*(p)$.
This is our main result.

\begin{theorem}\label{115}
  Let $p$ be an environment satisfying \reff{ccc} such that
  $\rho^*(p)<\infty$ and $\eta$ be a configuration such that
  $\underline D(\eta) > \rho^*(p)$. Then
\[
\lim_{t\to\infty} \delta_\eta S_p(t)=\nu_{p,c}\,,
\]
where $\delta_\eta$ is the measure giving weight one to the
configuration $\eta$.
\end{theorem}

\vskip 3mm As a corollary to Theorem \ref{115} we obtain the
asymptotic behavior of the system when the environment is randomly
chosen. Let $m$ be the distribution of a stationary ergodic
sequence $p$ on $[c,1]$ such that $m(\{p:p_0=c\})=0$, $m(\{p:
c<p_0 < c +\varepsilon \})>0$ for all $\varepsilon>0$. The measure
$m\nu_{\cdot,v}$ defined by $m\nu_{\cdot,v}f = \int m(dp)\int
\nu_{p,v}(d\eta) f(\eta)$ is an ergodic distribution on $\X$ and,
by the Ergodic Theorem, for all $v< c$ and for $m$-almost all $p$,
the asymptotic density exists $\nu_{p,v}$ a.s. and is equal to:
\[
R(v)=\int \frac {v}{p_0 -v}m(dp).
\]

Let $\rho^* := \lim_{v\to c} R(v)$ and assume $\rho^*<\infty$. In
this case for $m$-almost all environment $p$ any invariant measure
for $L_p$ is dominated by $\nu_{p,c}$.  The following theorem
concerns the behavior of the process when the initial measure
concentrates on configurations with density strictly higher than
$\rho^*$.

\begin{theorem}\label{2.4}
  Let $m$ be the distribution of a stationary ergodic sequence
  $p=(p_x)_{x\in\Z}$ on $(c,1]$ such that $m(\{p: c< p_0 < c
  +\varepsilon \})>0$ for all $\varepsilon>0$ and for which
  $\rho^*<\infty $. Let $\nu$ be a measure for which $\nu \ a.s.$
  $\underline D(\eta)$ is strictly bigger than $\rho^*$. Then, for
  $m$-almost all $p$
  \[
  \lim_{t\to\infty} \nu S_p(t)=\nu_{p,c}.
  \]
\end{theorem}

\section{Domination and Invariant measures}

In this section we prove Proposition \ref{2.2} and Theorem
\ref{2.1}.

\noindent{\bf Proof of Proposition 2.2.} Fix an arbitrary site $y$
and let $x_n$ be a decreasing sequence such that $x_1<y$, $p_{x_n}
<p_z$ for $x_n<z\le y$ and $p_{x_n}$ decreases to $c$. The
existence of such a sequence is guaranteed by \reff{ccc}. Consider
a process on ${\N}^{\{x_n +1,...,y\}}$ with generator given by:
\begin{eqnarray}
  \label{ln}
  L_{p,n}f(\eta )&=& \sum_{z=x_{n}+1}^{y-1}
{\bf 1 }\{\eta (z)>0\} p_z[f(\eta^z)-f(\eta)] \nonumber \\ &&\quad
+\; p_{x_n}[f(\eta+\mathfrak d_{x_n +1}) - f(\eta)] \\ &&\quad +\;
{\bf 1 }\{\eta (y)>0\} p_y [f(\eta -\mathfrak d_y)-f(\eta)] \;
.\nonumber
\end{eqnarray}
Let $S_{p,n}$ be the semigroup associated to this process and for
an arbitrary probability measure $\nu$ let $\nu _n$ be its
projection on ${\N}^{\{x_n +1,...,y\}}$. Standard coupling
arguments show that
\[
(\nu S_p(t))_n \leq \nu_n S_{p,n}(t) \; .
\]
The coupling of the two processes is done using the same Poisson
processes $N_x(t)$ defined in Section 2. The reason why the
domination holds is that for the process $S_{p,n}(t)$, each time
the Poisson process $N_{x_n}(t)$ jumps, a new particle appears in
$x_n+1$, while the same happens for the process $S_p(t)$ only when
there is at least a particle in the site $x_n$.

The process with generator $L_{p,n} $ is irreducible and  has a
countable state space, moreover a simple computation shows that
the product measure $\mu_{n,p}$ with marginals given by
\[
\mu _{n,p} \{\eta :\eta (z)=k\} = \Bigl(1-\frac
{p_{x_n}}{p_z}\Bigr)\Bigl(\frac {p_{x_n}}{p_z}\Bigr)^k,
\]
where $x_n<z\leq y $, is invariant for the process. Therefore
$\nu_n S_{n,p}(t)$ converges to $\mu _{n,p}$ and any weak limit
point of $(\nu S_p(t))_n$ is bounded above by $\mu_{n,p}$. Since
as $n$ goes to infinity the marginals of $\mu_{n,p}$ converge to
the marginals of $\nu_{c,p}$ the proposition is proved. $\carn$

\bigskip
\noindent{\bf Proof of Theorem \ref{2.1}.} Since only the final
step of the proof is different from the proof of Theorem 1.11 in
Andjel (1982) (in which the set of all invariant measures is
characterized for a family of asymmetric zero-range process) we
refer the reader to that paper. Exactly as there one proves that
if $\nu _p$ is an extremal invariant measure then for each
$v<\inf_x p_x$ either $\nu _p\leq \nu _{v,p}$ or $\nu _p\geq \nu
_{v,p}$. This implies that either $\nu _p =\nu _{v,p}$ for some
$v$ or $\nu _p \geq \nu _{v,p}$ for all $v$.  The latter case
cannot occur if there exists $x$ such that $p_x=\inf_y p_y$
because this would imply that $\nu _p \{\eta\ :\eta (x)>k \}=1$
for all $k$.  Therefore $\nu_p \geq \nu_{c,p}$ and either
$\liminf_{x\to -\infty}\ p_x=\inf p_y$ or $\liminf_{x\to \infty}\
p_x=\inf p_y$.  In the first of these cases, Proposition \ref{2.2}
allows us to conclude immediately. In the second case we argue by
contradiction: let $\widetilde\nu$ be a probability measure on
${\N}^{\Z}\times{\N}^{\Z}$ admitting as first marginal and second
marginal $\nu _p$ and $\nu _{c,p}$ respectively and such that
$\widetilde\nu \{(\eta ,\xi): \eta \geq \xi \}=1$. Consider the
standard coupled process with initial measure $\widetilde\nu$.
Denote by $\overline{S}(t)$ the semigroup associated to this
process and assume that for some $x$, $\widetilde\nu \{(\eta
,\xi): \eta (x) > \xi (x) \}>0$. Suppose it exists $k$ and $l$ in
$\N\setminus\{0\}$ such that
$\nu\{\eta(x)=k+l,\xi(x)=k\}=\varepsilon_1>0$ then at any time
$\delta>0$ one can find a $\varepsilon_2>0$ such that
\[
\nu\bar{S}(\delta)\{\eta(x)=k+l-1,\xi(x)=k-1\}=\nu\{\eta(x)=k+l-1,\xi(x)=k-1\}
=\varepsilon_2.
\]
To see that, one has just to control the arrivals and departures
of particles on sites $x-1$ and $x$ which are given by exponential
clocks. By induction it follows that for all $t>0$
\[
\widetilde\nu \overline{S}(t) \left\{(\eta ,\xi): \eta (x) > \xi
(x)=0 \right\}>0.
\]
Hence
\[
\nu _p\left\{ \eta :\eta (x)>0\right\}> \nu _{c,p}\left\{ \eta
:\eta (x)>0\right\}=\frac {c}{p_x}.
\]
Pick $y>x$ and such that $p_y < p_x\nu _p\{ \eta :\eta (x)>0\}$.
Then let $f(\eta )=\sum_{z=x+1}^y \eta (z)$. Now  a simple
calculation shows that $\int L_p f(\eta) d\nu_p(\eta)>0$
contradicting the invariance of $\nu_p$.$\carn$

\medskip
\noindent {\bf Remark:} Proofs of Theorem \ref{2.1} and
Proposition \ref{2.2} can easily be extended to a larger class of
one-dimensional nearest-neighbors asymmetric zero range processes
in non-homogeneous environment. In these systems a particle at
site $x$ on configuration $\eta$ jumps at rate $p_xg(\eta(x))$ to
site $x+1$, where $g:\N\to[0,\infty)$ is a non-decreasing bounded
function such that $g(0)=0$.

\section{Convergence}

We prove in this section Theorem \ref{115}. Fix a measure $\nu$ on
$\X$ concentrated on configurations with lower asymptotic left
density strictly greater than $\rho^*(p)$. Let $\tilde\nu$ be a
weak limit of $\nu S_p(t)$. Proposition \ref{2.2} shows that
$\tilde\nu$ is dominated by $\nu_{p,c}$. Lemma \ref{eq:lim} below
implies that $\tilde\nu$ dominates $\nu_{p,v}$ for all $v<c$. This
finishes the proof of Theorem \ref{115} because $\{\nu_{p,v}:0\le
v<c\}$ is an increasing sequence of measures converging to
$\nu_{p,c}$.

Denote $\{\overline{S}_p(t):\, t\ge 0\}$ the semigroup
corresponding to the coupling between two versions of the process
with (possibly) different initial configurations, by using the
same Poisson processes $(N_x(t):x\in\Z)$ in its construction.

\begin{lemma}
\label{eq:lim} Let $p$ be an environment satisfying \reff{ccc} and
such that $\rho^*(p)<\infty$ and $\zeta$ a configuration with
lower asymptotic left density $\underline D(\zeta)>\rho^*(p)$.
Then for any $v<c$,
\begin{equation}
  \label{4.1}
  \lim_{t\to\infty} (\delta_\zeta\times\nu_{p,v})\overline{S}_p(t)
\Big \{(\eta,\xi):\, \eta (x)<\xi(x) \Big \} \; =\; 0
\end{equation}
for all $x$ in ${\Z}$.
\end{lemma}

The proof of this lemma requires the following result. It states
that for each $v<c$ the asymptotic velocity of a second class
particle in the zero-range process in the environment $p$ under
the invariant measure $\nu_{p,v}$ is strictly positive.

Fix a starting site $z$ and consider a coupled zero range process
with initial condition $(\eta,\eta+\mathfrak d_{z})$ and semigroup
$\overline S_p(t)$. Under the coupled dynamics the number of sites
where the marginals differ does not increase in time. Let $X^z_t$
be the site where the marginals differ at time $t$. We can think
that $X^z_t$ stands for the position of a ``second class
particle''. Indeed, if the second class particle is at $x$ at time
$t$ it jumps to $x+1$ at rate $p_x\,\one\{\eta_t(x)=0\}$. In other
words, the second class particle jumps only if there is no other
particle at the site where it is.

For an environment $p$ and a probability measure $\nu$ on $\X$,
denote by $\P_\nu$ the measure on $D(\R_+, \X)$ induced by $\nu$
and the Markov process with generator $L_p$ defined in \reff{a3}.
In the next lemma we write $\P_{(\nu,z)}$ for a coupled process
whose initial configuration is $(\eta,\eta+\mathfrak d_z)$, with
$\eta$ distributed according to $\nu$. Since $R(p,v)$ is convex
and strictly increasing
\begin{equation}
  \label{602}
  \ga(p,v) := [R'(p,v)]^{-1}
\end{equation}
exists in a dense subset of $(0,c)$. In the sequel we abuse
notation by not writing integer parts where necessary.

\begin{lemma}
\label{secondclass} Let $p$ be an environment for which the limits
in \reff{cc1} exist for $v<c$. Pick $v\in (0,c)$ such that $
\ga(p,v)$ exists. Then,
\begin{equation}
\label{600} \lim_{t\to\infty}
\P_{(\nu_{p,v},-at)}\Bigl(\Bigl|\frac{X^{-at}_t}{t}-
(\ga(p,v)-a)\Bigr|\,>\vep\Bigr)\; = \;0\,,
\end{equation}
for all $\vep>0$ if $a> \ga(p,v)$ and
\begin{equation}
\label{601} \lim_{t\to\infty}
\P_{(\nu_{p,v},-at)}\Bigl(\frac{X^{-at}_t}{t} \ge 0\Bigr) \,=\, 1
\end{equation}
if  $a<\ga(p,v)$.
\end{lemma}

\vskip 3mm \noindent{\bf Remark.} The more complete result when
the starting point $a$ is greater than $\ga(p,v)$ comes from the
fact that in our hypothesis we have only the asymptotic
\emph{left} limits \reff{cc1}. If the limits \reff{cc1} hold for
both sides, then \reff{600} is valid for all $a$.

\noindent{\bf Proof:} Note that it suffices to prove (\ref{600}),
since (\ref{601}) follows from (\ref{600}) because $X_t^x\le
X_t^y$ for all $t\ge 0$ if $x\le y$ and because (\ref{601}) does
not depend on the environment to the right of the origin. For
$u<w<c$ let $\bar \nu_{p,u,w}$ be the product measure on $\X\times
{\X}$ whose first marginal is equal to $\nu_{p,u}$, whose second
marginal is equal to $\nu_{p,w}$ and which is concentrated above
the diagonal: $\bar \nu_{p,u,w} \{(\eta, \xi):\, \eta \le
\xi\}=1$.  Denote by $(\eta_t, \xi_t)$ the coupled Markov process
starting from $\bar \nu_{p,u,w}$.

Denote by $\zeta_t$ the difference $\xi_t - \eta_t$ and observe
that the $\zeta$-particles evolve as second class particles in the
sense that a $\zeta$-particle jumps from $x$ to $x+1$ at rate
\[
p_x\, [\one\{\eta (x) + \zeta (x)\ge 1\} - \one\{\eta (x)\ge 1\}];
\]
that is, when there are no $\eta$ particles present.  In this case
we say that the $\eta$ particles have \emph{priority} over the
$\zeta$ particles. We label the $\zeta$-particles at time $0$ in
the following way. Without losing much we can assume that there is
a $\zeta$ particle at site (integer part of) $-at$. The measure
conditioned on this event is absolutely continuous with respect to
$\bar \nu_{p,u,w}$, an this will be enough for our purposes, as we
shall only use laws of large numbers. Call particle $0$ this
particle, and complete the labeling in such a way that a particle
with label $j$ is at the same site or at the left of a particle
with label $k$ if $j<k$. Denote by $Y^j_t$ the position at time
$t$ of the particle labeled $j$. By construction, we have $\cdots
\le Y^{-1}_0 < Y^{0}_0 = -at \le Y^1_0 \le \cdots$. We let the
second class particles evolve in a way to preserve this order. To
keep track of the densities involved in the definition we call
$Y^{u,w}_t = Y^0_t$.

Consider now a single second class particle for the $\eta$ process
initially at the position of $-at$. This is obtained by
considering the coupled initial condition $(\eta,\eta+\mathfrak
d_{-at})$. Denote the position of the single second class particle
at time $t$ by $X^u_t$ (for $u=v$, this has the same law as the
particle denoted by $X^{-at}_t$ in the statement of the
proposition). Since $Y^{u,w}_0= X^u_0$, in the coupled evolution
obtained by using the same Poisson processes $(N_x(t))$ we have
$Y^{u,w}_t\le X^u_t$ for all $t$. Indeed, in this coupling $Y^k_t$
for $k>0$ have priority over $Y^{u,w}_t$ while those particles
have no priority over $X^u_t$. Similarly, consider a second class
particle for the $\xi$ process and denote it $X^w_t$.  Since
$Y^k_t$ for $k<0$ have priority over $X^w_t$ but not over
$Y^{u,w}_t$, $X^w_t\le Y^{u,w}_t$.  Hence, for $0\le u<w\le c$,
\begin{equation}
  \label{0vwc}
  X^w_t\le Y^{u,w}_t\le X^u_t\,,
\end{equation}
$\P_{\bar \nu_{p,u,w}}$ almost surely.

Denote by $J_t^1$, $J_t^{1+2}$ and $J_t^2$ the total number of
$\eta$, $\xi$ and $\zeta$ particles that jumped from $-at$ to
$-at+1$ before time $t$. In particular, $J_t^2 = J_t^{1+2} -
J_t^1$. By Burke's theorem, the number of $\eta$-particles (resp.
$\xi$-particles) that jump from $-at$ to $-at+1$ is a Poisson
process of parameter $u$ (resp. $w$). Hence the number of
$\zeta$-particles that jump from $-at$ to $-at+1$ in the interval
$[0,t]$ is the difference of two Poisson processes and satisfies
the law of large numbers:
\[
\lim_{t\to \infty} \frac{J_t^2}t \; =\; \lim_{t\to \infty}
\frac{J_t^{1+2} - J^1_t}t \; =\; w  - u
\]
in $\P_{\bar \nu_{p,u,w}}$ probability. On the other hand, for
every $t\ge 0$,
\[
J_t^2 \; := \; \sum_{x=-at+1}^{Y^{u,w}_t}\zeta_t(x)-A_t \; = \;
\sum_{x=-at+1}^{Y^{u,w}_t} \xi_t(x) \; - \;
\sum_{x=-at+1}^{Y^{u,w}_t} \eta_t(x)-A_t \; .
\]
where  $\vert A_t\vert \leq \zeta_t(Y^{u,w}_t)$. Note that
$\zeta_t( Y^{u,w}_t)$ is stochastically bounded above by a
geometric random variable of parameter $w/c$. Therefore $\vert
{A_t}/{t}\vert$ converges to $0$ in $\P_{\bar \nu_{p,u,w}}$
probability as $t$ goes to infinity. As in the proof of Theorem
12.1 of Ferrari (1992), it follows from the previous equation and
the law of large numbers for $\eta_t$ and $\xi_t$ ---that are
distributed according to product (invariant) measures with
densities $R(p,u)$ and $R(p,w)$ respectively--- that for $u$ and
$w$ strictly smaller than $c$,
\begin{equation}
  \label{ytp}
  \lim_{t\uparrow\infty} {Y^{u,w}_t\over t} \; +\; a
\;=\;{w-u \over R(p,w) - R(p,u)}
\end{equation}
in $\P_{\bar \nu_{p,u,w}}$ probability. Notice that we used here
the fact that $a>\gamma(p,v)$. In this case $Y_t^{u,w}/t<0$ and
the previous sums refer only to negative sites. Hence, from
\reff{0vwc} we have
\begin{equation}
  \label{xcw}
  \lim_{t\uparrow\infty} {X^{w}_t\over t} \; +\; a \;\le\;{w-u \over R(p,w) -
  R(p,u)} \;\le\;\lim_{t\uparrow\infty} {X^{u}_t\over t} \; +\; a
\end{equation}
in $\P_{\bar \nu_{p,u,w}}$ probability. Fixing $w=v$ and taking
the limit $u\to v$ and then fixing $u=v$ and taking the limit
$w\to v$ and taking account of the differentiability of $R$ in
$v$, we get \reff{600} and \reff{601}. $\carn$

\bigskip
We are now in a position to prove Lemma \ref{eq:lim}.

\noindent {\bf Proof of lemma \ref{eq:lim}} The proof is performed
via coupling. We start with two different initial configurations
$\eta$ and $\xi$ with marginal distributions $\nu$ and
$\nu_{p,v}$, respectively. Hence $\eta$ has lower asymptotic
density bigger than $\rho^*(p)$ and $\xi$ has asymptotic density
$R(p,v)$. We use the same Poisson processes for both processes and
call $(\eta_t,\xi_t)$ the coupled process. The configurations
$\eta$ and $\xi$ are in principle not ordered: there are (possibly
an infinite number of) sites $z$ such that $(\eta(z)-\xi(z))^+>0$
and (possibly an infinite number of) sites $y$ such that
$(\eta(y)-\xi(y))^->0$. We say that we have $\eta\xi$
discrepancies in the first case and $\xi\eta$ discrepancies in the
second. The number of coupled particles at site $x$ at time $t$ is
given by
\begin{equation}
  \label{cp}
 \bar\xi_t(x):= \min \{\eta_t(x),\xi_t(x)\}
\end{equation}
The $\bar\xi$ particles move as regular (first class) zero range
particles. There is at most one type of discrepancy at each site
at time zero. Discrepancies of both types move as second class
particles with respect to the already coupled particles. When a
$\eta\xi$ discrepancy jumps to a site $z$ occupied by at least one
$\xi\eta$ discrepancy, the $\eta\xi$ discrepancy and one of the
$\xi\eta$ discrepancies at $z$ coalesce into a coupled $\bar\xi$
particle in $z$. The coupled particle behaves from this moment on
as a regular (first class) particle. The same is true when the
roles of $\xi$ and $\eta$ are reversed.

The above description of the evolution implies in particular that
a tagged discrepancy can not go through a region occupied by the
other type of discrepancies.

We will choose a negative site $y$ such that the jump rate from
$y-1$ to $y$ is close to $c$. Then we follow the $\xi\eta$
discrepancies belonging to two disjoint regions of ${\Z}$ at time
$0$ and give upper bounds on the probability of finding them at
$y$ at time $t$.

Roughly speaking, a $\xi\eta$ discrepancy at $y$ cannot come from
a region ``close'' to $y$ because we prove that there is a minimum
positive velocity for the $\xi\eta$ discrepancies to go. This
velocity is given by the velocity of a second class particle under
$\nu_{p,v}$. On the other hand, the $\xi\eta$ discrepancy cannot
come from a region ``far'' from $y$ because due to the difference
of densities, a lot of $\eta\xi$ discrepancies will be between it
and $y$ and hence they must pass site $y-1$ before it. But since
we have chosen a small rate for this site, a traffic rush will
prevent them to pass.  With this idea in mind, we have to choose
the ``close'' and ``far'' regions and the value of the rate
at~$y-1$.

Fix $v<c$ such that $R(p,\cdot)$ is differentiable in $v$. Let
$\gamma= \gamma(p,v)$, the (strictly positive) asymptotic speed of
a second class particle under $\nu_{p,v}$ in the sense
of~(\ref{602}). Denote by $\beta$ the difference between the lower
asymptotic density of $\eta$ and $R(p,v)$~:
\[
\beta=\beta(p,v) \; =\; \liminf_{n\to \infty} \frac 1n
\sum_{x=-n+1}^0 [\eta (x) - R(p,v)]\; .
\]
For reasons that will become clear later (cf.\/ display
\reff{lnp}), we let
\[
b\,=\,b(p,v)\,=\, {R'(p,v) (c-v) \over \rho^*(p)-R(p,v)} \,<\,1\,,
\]
by the convexity of $R$; recall that $\rho^*(p) = \lim_{v\to c}
R(p,v)$. With this choice,
\begin{eqnarray}
  \label{gb4}
   \beta\gamma b - c + v \;=\; { c-v \over \rho^*(p)-R(p,v) } \Big\{
  \beta - [\rho^*(p) - R(p,v)]\Big\} \; >\; 0\; .
\end{eqnarray}
This allows us to fix $\varepsilon = \varepsilon(v)$ satisfying
\[
0 \,<\, \varepsilon(v)\, <\, \beta \gamma b - c + v\,.
\]

Finally, choose a negative site $y=y(v)$ such that
\begin{equation}
  \label{y11}
  p_{y-1}< c+ \varepsilon\,.
\end{equation}
We shall prove that
\begin{equation}
\label{a1} \lim_{t\to\infty}\; (\nu\times\nu_{p,v})\overline{S}_p
(t) \Big \{(\eta,\xi):\, \eta (y)<\xi(y) \Big \} \; =\; 0 \; .
\end{equation}

We can order the $\xi\eta$ discrepancies and assume without loss
of generality that the order is preserved in future times as we
did in Lemma \ref{secondclass}. Of course some of the
discrepancies will disappear. Let $Z^k=Z^k_t(\xi,\eta)$ the
positions of the ordered $\xi\eta$ discrepancies at time $t$ with
the convention that $Z^k_t=\infty$ if the corresponding
discrepancy coalesced with a $\eta\xi$ one giving place to a
$\bar\xi$ coupled particle. Let
\begin{eqnarray}
  \label{agt}
 \lefteqn{ A_{\ga,t}(\eta,\xi)}\nonumber \\
&:=& \left\{\hbox{a $\xi\eta$ discrepancy in the box $[y-(t\gamma
\bar b),y]$ at time $0$ }\right.\\ &&\ \left.\hbox{has moved to
site $y$ at time $t$}\right\} \nonumber\\ &:=&
\cup_k\left\{Z^k_0\in\left[y-(t\gamma \bar b),y\right],
\;Z^k_t=y\right\} \nonumber
\end{eqnarray}
where $\bar b\,:=\,(1+b)/2\,\in\,(b,1)$. Hence
\begin{equation}
  \label{lzp}
  \P(A_{\ga,t}(\eta,\xi)) \,\le\, \P\left(\min\left\{Z^k_t: Z^k_0\in
  \left[y-(t\gamma \bar b),y\right]\right\} \le \,y\right)
\end{equation}

We wish to give and upper bound to the event in the right hand
side above. To do so we consider the coupled $(\eta , \xi) $
process and the $\xi $ process to which we add a unique second
class particle at $y-(t\gamma \bar b)$, evolving together with
jumps occurring at times given by the same Poisson processes. We
denote by $X_t^{y-t\gamma\bar b}$ the position of the second class
particle at time $t$.  If the second class particle has reached $
y+1$ no later than time $t$, then there exists an increasing
sequence of random times $0<T_{y-(t\gamma \bar b
  )}<T_{y-(t\gamma \bar b)+1}<...<T_y$ such that at each of these
times the corresponding site has been emptied of its $\xi $
particles. But this implies that all the $\xi \eta $ discrepancies
which at time $0$ were in the interval $[y-(t\gamma \bar b),y]$
have disappeared or are strictly to the right of $y$. Therefore:
\begin{equation}
  \label{aia}
  \P\left(A_{\ga,t}(\eta,\xi)\right)\,\le\,
  \P(X^{y-t\gamma \bar b}_t \le y )\,.
\end{equation}
By \reff{601} this tends to $0$ as $t$ tends to infinity because
$\gamma \bar b<\ga$. The above argument is independent of the
value of $p_{y-1}$.

It now suffices to check that the probability that a $\xi\eta$
discrepancy, to the left of $y-(t\gamma\bar b)$ at time 0 reaches
$y$ no later than time $t$, tends to $0$ as $t$ tends to infinity.
Let
\begin{eqnarray}
  \label{bgt}
 \lefteqn{ B_{\ga,t}(\eta,\xi)}\nonumber \\
&:=& \left\{\hbox{a $\xi\eta$ discrepancy in $(-\infty,y-(t\gamma
\bar b)]$ at time $0$}\right.\\ &&\ \left.\hbox{ has moved to site
$y$ at time $t$}\right\} \nonumber\\ &:=&
\cup_k\left\{Z^k_0\in(-\infty,y-(t\gamma \bar b)],
\;Z^k_t=y\right\}\nonumber
\end{eqnarray}
Call $W_t^k(\eta,\zeta)$ the positions of the $\eta\xi$
discrepancies at time $t$, $W_0^0$ being the first $\eta\xi$
discrepancy to the left of the origin. As before set
$W^k_t=\infty$ if the $k$th discrepancy coalesced with a $\xi\eta$
one before $t$.

Since a  $\xi \eta $ discrepancy cannot cross over an $\eta \xi $
discrepancy,
\begin{eqnarray}
\nonumber \lefteqn{\hskip-2cm {B_{\gamma,t}(\eta,\xi)\cap
  \left( \bigcap_{z\le y-\gamma t\bar b} \left\{\sum_{x=z}^{y-1}
    \left(\eta_0(x)-\xi_0(x)\right) > \gbt b\right\} \right)}}\\
  &&\qquad\qquad\qquad
\subset \left\{ I^2_t - I^1_t >\gbt b\right\}  \label{weh}
\end{eqnarray}
where $I^2_t$ and $I^1_t$ are the number of $\eta$, respectively
$\xi$, particles jumping from $y-1$ to $y$ in the interval
$[0,t]$. Since
\begin{eqnarray}
  \label{100}
 \lefteqn{ \left\{ \bigcap_{z\le y-\gamma t\bar b}
\left\{\sum_{x=z}^{y-1}
    (\eta_0(x)-\xi_0(x)) > \gbt b\right\}  \right\}^c}\nonumber\\
&=& \bigcup_{z\le y-\gamma t\bar b} \left\{\sum_{x=z}^{y-1}
    \left(\eta_0(x)-\xi_0(x)\right) \leq  \gbt b\right\},
\end{eqnarray}

\noindent to bound $\P(B_{\gamma,t}(\eta,\xi))$ it suffices to
bound the probabilities of the sets on the right hand sides of
\reff{weh} and \reff{100}. For \reff{weh} we have
\begin{equation}
  \label{101}
  \P(I^2_t - I^1_t >\gbt b)\,\le\,\P(N^{c+\varepsilon}_t \,-
  \,N^v_t>\gbt b),
\end{equation}
where $N^a_t$ is a Poisson process of parameter $a$. The above
inequality holds because the $\eta$-particles jump from $y-1$ to
$y$ at rate not greater than $p_{y-1}$, which is by construction
less than or equal to $c+\varepsilon$. On the other hand, by
Burke's theorem, the number of jumps from $y-1$ to $y$ for the
$\xi$-particles is a Poisson process of rate $v$. By the law of
large numbers for the Poisson processes, we have
\begin{eqnarray}
  \label{lnp}
  \lim_{t\to\infty} {1\over t} (N^{c+\varepsilon}_t \,-
  \,N^v_t) &=& c-v+\varepsilon \;<\;\beta\gamma b\,,
\end{eqnarray}
because we chose $\varepsilon < \gamma\beta b -c +v$. Hence
\reff{101} goes to zero as $t\to\infty$.

On the other hand, the probability of the set in the right hand
side of \reff{100} is
\begin{equation}
  \label{102}
 \P \left(\sup_{z \leq y-\gamma t\bar b} \sum_{x=z}^{y-1}
(\eta_0(x)-\xi_0(x)) \leq  \gbt b\right)\end{equation} By the
ergodicity of $\xi$ and the fact that $\eta$ has left density,
with probability one:
\begin{eqnarray}
\lim_{t\to\infty}\;{1\over t} \sum_{x=y-t\gamma \bar b}^{y-1}
(\eta_0(x)-\xi_0(x))&=& {\gamma \beta \bar b} \;>\; \gamma\beta
b\,,
\end{eqnarray}
by the way we chose $\bar b$. This implies that \reff{102} goes to
zero as $t\to\infty$. This proves (\ref{a1}).

To deduce the statement of the lemma from (\ref{a1}) we need the
following lemma which says that if there exists a subsequence of
times giving positive probability to a cylinder set, then any
other cylinder set obtained by moving one particle to the right
has the same property. These lines follow Andjel~(1982).

\begin{lemma}
  \label{andj}
Let $f$ be the following cylinder function on $\N^\Z\times\N^\Z$.
\begin{equation}
  \label{201}
  f(\eta,\xi) = \one\{\eta(x) = \bar\eta(x),\, \xi(x) = \bar\xi(x)\,: x\in A\}
\end{equation}
for some finite $A\subset \Z$ and arbitrary configurations
$\bar\eta,\bar\xi \in \N^\Z$. Let $z\in\Z$ be an arbitrary site;
define $f^z$ as
\begin{equation}
  \label{2011}
  f^z(\eta,\xi) = \one\{\eta(x) = \bar\eta^z(x),\, \xi(x) = \bar\xi^z(x)\,:
  x\in A\}
\end{equation}
Let $(\eta_t,\xi_t)$ be the coupled process starting from an
arbitrary measure. Then
\begin{equation}
  \label{andj1}
  \limsup_{t\to\infty} \E f(\eta_t,\xi_t) > 0 \;\;\hbox{ implies } \;\;
  \limsup_{t\to\infty}\E f^z(\eta_t,\xi_t) > 0
\end{equation}
\end{lemma}

\proof Let $\tilde A = \{x\in\Z\,:\, x+1\in A\}$. Since if
$z\notin A\cup\tilde A$ implies that $f(\eta,\xi)=f^z(\eta,\xi)$
(and hence for these $z$ the lemma is trivial), we fix a $z\in
A\cup\tilde A$. Assume that $t_n$ is a sequence of times such that
\begin{equation}
  \label{202}
  \lim_{n\to\infty} \E f(\eta_{t_n},\xi_{t_n}) = c>0
\end{equation}
Fix a time $s$ (equal to one, for instance) and consider the event
$\B_n = \{N_z(t_n+s)-N_z(t_n) = 1,\, N_x(t_n+s)-N_x(t_n)=0
\;\hbox{ for }x\in A\cup\tilde A\setminus \{z\}\}$. That is the
event ``exactly one Poisson event occurs for $z$ in the interval
$[t_n, t_{n}+s)$ and no events occur for the other sites in
$A\cup\tilde A$ in the same time interval''. Then
\begin{eqnarray}
  \label{203}
 && \E f(\eta_{t_n},\xi_{t_n})\,\P(\B_n)\;
 \le\; \E f^z(\eta_{t_n+s},\xi_{t_n+s})
\end{eqnarray}
Since the probability of $\B_n$ is independent of $n$ and
positive, this proves the lemma. \quad $\carn$

We continue with the proof of Lemma \ref{eq:lim}. Take an
arbitrary $y$ satisfying \reff{a1}.  Consider the coupled process
starting with the measure $(\nu\times\nu_{p,v})$. By Proposition
\ref{2.2} both $(\eta_t(y-1),\xi_t(y-1))$ and
$(\eta_t(y),\xi_t(y))$ are tight sequences.  Hence there exists a
$K$ depending on $p_{y-1}, p_y$ such that
\begin{equation}
  \label{205}
 \lim_{n\to\infty} \P(\eta_{t_n}(y-1)<\xi_{t_n}(y-1))>0
\end{equation}
implies
\begin{equation}
  \label{206}
\lim_{n\to\infty}\P(\eta_{t_n}(y-1)<\xi_{t_n}(y-1),\;
\xi_{t_n}(y-1)\le K,\; \eta_{t_n}(y)\le K)>0
\end{equation}
Now we apply Lemma \ref{andj}: move first the (at most) $K$ $\eta$
particles from $y$, then the (at most) $K$ extra $\xi$ particles
from $y-1$ to $y$ to obtain that \reff{206} implies
\begin{equation}
  \label{207}
 \lim_{n\to\infty} \P(\eta_{t'_n}(y)<\xi_{t'_n}(y))>0
\end{equation}
for some subsequence $(t'_n)$, in contradiction with~(\ref{a1}).
With the same argument we can go to $x=y-2, y-3, \dots$. This
proves that \reff{4.1} holds for all $x<y$ for $y$ satisfying
\reff{y11}. On the other hand, the marginal law of the coupled
process at $x$ does not depend on the value of $p_y$ for $y>x$.
Hence, we can assume \reff{y11} for $y\ge x+2$ and obtain the
result for all $x\in\Z$. This argument works because when we
modify $p_{y-1}$ we change the process only to the right of $y-1$,
maintaining the values $R(p,v)$ and $\gamma(p,v)$ unaltered, as
they are asymptotic \emph{left} values. For this reason we can use
the same $\varepsilon$ in \reff{y11}.  \quad $\carn$

\bigskip

{\bf Acknoledgements:} The authors thank the referee for his/her
cautious reading. PAF and HG would like to thank Joachim Krug for
fruitful discussions. PAF and HG thank FAPESP, PROBAL/CAPES and
FINEP-PRONEX for their support.

\end{document}